\documentclass{amsart}
\usepackage{amssymb}

\def\TR{\operatorname{Tr}}
\newtheorem{Theorem}{Theorem}

\title{Enumerating Permutation Polynomials over finite fields
by degree}

\author{Sergei Konyagin}

\address[Konyagin]{Department of Mechanics and Mathematics, Moscow State
University, Vorobjovy Gory, 119899 Moscow, Russia}
\email[Konyagin]{kon@mech.math.msu.su, ars204@arstel.ru}

\author{Francesco Pappalardi}
 \address[Pappalardi]{Dipartimento di Matematica,
Universit\`{a} degli Studi Roma Tre, Largo S. L. Murialdo, 1,
I--00146 Roma, Italia}
\email[Pappalardi]{pappa@mat.uniroma3.it}

\begin{document}

\begin{abstract}
We prove an asymptotic formula for the number of permutation for
which the associated permutation polynomial has degree smaller
than $q-2$.
\end{abstract}
\maketitle

Let $\mathbb F_q$ be a finite field with $q=p^f>2$ elements and
let $\sigma\in\mathcal S(\mathbb F_q)$ be a permutation of the
elements of $\mathbb F_q$. The \text{permutation polynomial}
$f_\sigma$ of $\sigma$ is
$$f_\sigma(x)=\sum_{c\in \mathbb F_q} \sigma(c)\left(1-(x-c)^{q-1}
\right)\in\mathbb F_q[x].$$
$f_\sigma$ has the  property that $f_\sigma(a)=\sigma(a)$ for every $a\in
\mathbb F_q$ and this explains its name.

For an account of the basic properties of permutation polynomials
we refer to the book of Lidl and  Niederreiter \cite{LN}.

{From} the definition, it follows that for every $\sigma$
$$\partial(f_\sigma)\leq q-2.$$

A variety of problems and questions regarding Permutation
polynomials have been posed by Lidl and Mullen in \cite{LM1,LM2}.
Among these there is problem of determining the number $N_d$ of
permutation polynomials of fixed degree $d$. In \cite{CP} and
\cite{CP1}, Malvenuto and the second author address the problem of
counting the permutations that move a fixed number of elements of
$\mathbb F_q$ and whose permutation polynomials have ``low''
degree.

Here we consider all permutations and we want to prove the following

\begin{Theorem}
Let
$$N=\#\left\{\sigma\in\mathcal S(\mathbb F_q)\ |\
\partial(f_\sigma)< q-2\right\}.$$
Then
$$|N-(q-1)!| \leq \sqrt{2e/\pi}q^{q/2}.$$
\end{Theorem}

This confirms the common believe that \textsl{almost all
permutation polynomials have degree $q-2$}.
\bigskip

The first few values of $N$ are listed below:\bigskip

\begin{center}
\begin{tabular}{|r|r|r|r|r|r|r|c|c|} \hline
      $q$ & 2 & 3 &  4 &  5 &   7 &    8 &   9   &   11    \\
\hline
      $N$ & 0 & 0 & 12 & 20 & 630 & 5368 & 42120 & 3634950 \\
\hline
 $(q-1)!$ & 1 & 2 &  6 & 24 & 720 & 5040 & 40320 & 3628800 \\
\hline \end{tabular} \end{center}\bigskip

\noindent\textbf{Proof.} The proof uses exponential sums and
a similar argument as the one in \cite{K}.

By extracting the coefficient of $x^{q-2}$ in $f_\sigma(x)$, we
obtain that the degree of $f_\sigma(x)$ is strictly smaller than
$q-2$ if and only if
$$\sum_{c\in \mathbb F_q} c\sigma(c)=0.$$

For a fixed subset $S$ of $\mathbb F_q$, we introduce the
auxiliary set of functions
$$N_S=\left\{f\ |\ f: \mathbb F_q\longrightarrow S,
\textrm{ and } \sum_{c\in  S} cf(c)=0 \right\}$$ and set $n_S=\#
N_S$. By inclusion exclusion, it is easy to check that
\begin{equation}\label{uno}
N=\sum_{S\subseteq \mathbb F_q}(-1)^{q-|S|}n_S.
\end{equation}

Now if $e_p(u)=e^{\frac{2\pi i u}{p}}$, consider the identity
$$n_S=\frac{1}{q}\sum_{a\in\mathbb F_q}\left(
\sum_{f: \mathbb F_q\longrightarrow S}
e_p(\sum_{c\in\mathbb F_q}\TR(acf(c)))
\right)$$
which follows from the standard property
$$
\frac{1}{q}\sum_{a\in\mathbb F_q}
e_p(\TR(ax))=\begin{cases} 1 & \textrm{if } x=0\\ 0, &
\textrm{if } x\neq0. \end{cases}
$$

By exchanging the sum with the product, we obtain
$$n_S=\frac{1}{q}\sum_{a\in\mathbb F_q}\left(
\prod_{c\in \mathbb F_q}\sum_{t\in S}e_p(\TR(act)) \right).$$

By isolating the term with $a=0$ in the external sum and noticing
that the internal product does not depend on $a$ (for $a\neq0$),
we get
$$n_S=\frac{|S|^q}{q}+\frac{1}{q}
\sum_{a\in\mathbb F_q^*}\left(
\prod_{b\in \mathbb F_q}\sum_{t\in S}
e_p(\TR(bt))
\right).$$

Finally

\begin{equation}\label{due}
n_S=\frac{|S|^q}{q}+\frac{q-1}{q}
\prod_{b\in \mathbb F_q}\sum_{t\in S}
e_p(\TR(bt)).
\end{equation}

Now let us plug equation (\ref{due}) in equation (\ref{uno})
and obtain
$$
N-\sum_{S\subseteq \mathbb F_q}\frac{(-1)^{q-|S|}}{q}|S|^q =
\frac{q-1}{q}\sum_{S\subseteq \mathbb F_q}(-1)^{q-|S|}\prod_{b\in
\mathbb F_q}\sum_{t\in S} e_p(\TR(bt)).
$$

Note that inclusion exclusion gives
$$\sum_{S\subseteq \mathbb F_q}\frac{(-1)^{q-|S|}}{q}|S|^q=(q-1)!.$$

Therefore
$$
N-(q-1)! =\frac{q-1}{q}\sum_{S\subseteq \mathbb F_q}
(-1)^{q-|S|}|S|\prod_{b\in \mathbb F_q^*}\sum_{t\in S}
e_p(\TR(bt)).
$$

Using the fact that for $b \in\mathbb F_q^*$
$$\sum_{t \in S} e_p(\TR(bt)) =
-\sum_{t \not \in S}  e_p(\TR(bt))$$
and grouping together the term relative to $S$ and the term relative
to $\mathbb F_q\setminus S$, we get

\begin{equation}\label{bravoigor} \left| N-(q-1)! \right|\leq \frac{q-1}{2q}\sum_{S\subseteq
\mathbb F_q}\left|q-2|S|\right| \prod_{b\in \mathbb
F_q^*}\left|\sum_{t\in S} e_p(\TR(bt))\right|.
\end{equation}

Now let us also observe that
$$
\sum_{b\in \mathbb F_q}\left|\sum_{t\in S}
e_p(\TR(bt))\right|^2=q|S|,
$$
so that
$$
\sum_{b\in \mathbb F_q^*}\left|\sum_{t\in S}
e_p(\TR(bt))\right|^2=(q-|S|)|S|.
$$

{From} the fact that the geometric mean is always bounded by the
arithmetic mean (i.e. $(\prod_{i=1}^k |a_i|^2)^{1/k} \leq
\frac{1}{k}\sum_{i=1}^k |a_i|^2$), we have that
\begin{equation}\begin{array}{rcl}
\displaystyle{\prod_{b\in \mathbb F_q^*}\left|\sum_{t\in S}
e_p(\TR(bt))\right|}&\displaystyle{\leq}&
\displaystyle{\left(\frac{1}{q-1} \sum_{b\in \mathbb
F_q^*}\left|\sum_{t\in S} e_p(\TR(bt))\right|^2
\right)^{(q-1)/2}}\\
\\
&=&\displaystyle{\left(\frac{(q-|S|)|S|}{q-1}\right)^{(q-1)/2}}.
\end{array}\label{bravobis}\end{equation}

Furthermore, using (\ref{bravoigor}) and (\ref{bravobis}) we
obtain
\begin{equation}\left| N-(q-1)! \right|\leq
\frac{q-1}{2q(q-1)^{(q-1)/2}}\sum_{S\subseteq \mathbb
F_q}\left|q-2|S|\right| \left((q-|S|)|S|\right)^{(q-1)/2}.
\label{fine}\end{equation}

We want to estimate the above sum. Consider the inequality
\begin{equation}\label{finebis}
\left((q-|S|)|S|\right)^{(q-1)/2}\le
\left(\frac{q}{2}\right)^{q-1},\end{equation}
and the identity
\begin{equation}\label{bino}
\sum_{S\subseteq \mathbb F_q}\left|q-2|S|\right|
=2q\binom{q-1}{[q/2]},
\end{equation}
which holds since 
$$
\begin{array}{rcl}
2\sum_{\genfrac{}{}{0pt}{}{S\subseteq \mathbb F_q,}{|S|\leq q/2}}\left(q-2|S|\right)&=&
2\left[\sum_{j=0}^{[q/2]}\binom{q}{j}(q-j)-
\sum_{j=1}^{[q/2]}\binom{q}{j}(j)\right]\\
&=&2q \left[\sum_{j=0}^{[q/2]}\binom{q-1}{j}-
\sum_{j=1}^{[q/2]}\binom{q-1}{j-1}\right]=2q\binom{q-1}{[q/2]}.
\end{array}
$$

{From} the standard inequality
$$\binom{2n}{n}\leq\sqrt{\frac{2}{\pi}}\frac{2^{2n}}{\sqrt{2n+1/2}}$$
which can be found for example in \cite{GR}, we deduce
\begin{equation}\label{summ}
\binom{q-1}{[q/2]}\leq\sqrt{\frac{2}{\pi}}\frac{2^{q-1}}{\sqrt{q-1/2}}.
\end{equation}

Therefore, (\ref{fine}), (\ref{finebis}), (\ref{bino}) 
and (\ref{summ}) imply 
$$|N-(q-1)!|\leq \left(\frac{q-1}{\sqrt{q-1/2}\sqrt{q}}\right)
\sqrt{\frac{2}{\pi}}\left(\frac
q{q-1}\right)^{(q-1)/2}q^{q/2}$$
and in view of the inequalities
$$\frac{q-1}{\sqrt{q-1/2}\sqrt{q}}<1,\hspace{1cm}\left(\frac
q{q-1}\right)^{(q-1)/2}<\sqrt{e}$$ 
we finally obtain
$$\left|
N-(q-1)! \right|\leq
\sqrt{\frac{2e}{\pi}}{q^{q/2}}
$$
and this completes the proof.\hfill$_\blacksquare$\bigskip
\bigskip

\noindent\textbf{Conclusion.} Computations suggest that a more
careful estimate of the sum in (\ref{fine}) would yield to a
constant $\sqrt{\frac{e}{2\pi}}$ instead of $\sqrt{\frac{2e}{\pi}}$ as coefficient in $q^{q/2}$ in the statement of Theorem~1. However we feel 
that such a minor improvement does not justify the extra work.

The ideas in the proof of Theorem~1 can be used to deal with the
analogous problem of enumerating the permutation polynomials that
have the $i$--th coefficient equal to $0$ and also to the problem
of enumerating the permutation polynomials with degree less than
$q-k$ (for fixed $k$). However, the exponential sums that need to
be considered are significantly more complicated.\bigskip\bigskip

\centerline{\textbf{Acknowledgements.}}

The first author was supported by Grants 99-01-00357
and 00-15-96109 from the Russian Foundation for Basic Research.
The second author was partially supported by G.N.S.A.G.A. from
Istituto Nazionale di Alta Matematica.\medskip

The authors would like to thank Igor Shparlinski
for suggesting a substantial improvement with respect to the original
result.
\bigskip\bigskip

\end{document}